\documentclass[12pt]{article}

\usepackage[cp1251]{inputenc}
\usepackage[english]{babel}

\usepackage{amsfonts,epsfig,epstopdf}
\usepackage{amssymb,amsmath,enumerate,float}
\makeatother

\addto\captionsenglish{}

\newcommand{\ptl}{\partial}

\newcommand{\be}{\begin{equation}}
\newcommand{\ee}{\end{equation}}
\newcommand{\beq}{\begin{equation}}
\newcommand{\eeq}{\end{equation}}

\title{First arriving signals  in layered waveguides.
An approach based on dispersion diagrams}

\author{A. V. Shanin\\
M.V.Lomonosov Moscow State University, department of physics \\
119992, Russia, Moscow, Leninskie gory, MSU}

\begin{document}

\maketitle

\begin{abstract}

The first arriving signal (FAS) in a layered waveguide is investigated. It is well known
that the velocity of such a signal is close to the velocity of the fastest medium in the waveguide,
and it may be bigger than the fastest group velocity given by the dispersion diagram of the waveguide.
Usually the FAS pulse decays with the propagation distance. A model layered waveguide is studied
in the paper. It is shown that the FAS is associated with the pseudo-branch structure of the
dispersion diagram. The velocity of FAS is determined by the slope of the pseudo-branch. The decay is exponential
and it depends on the structure of the pseudo-branch.

\end{abstract}


\section{Introduction}

Consider a layered 2D acoustic waveguide, which is homogenous in the $x$-direction and has a sandwich structure
in the $y$-direction. The media constituting the waveguide can be elastic, liquid or gaseous. Waves in such a waveguide
are described by a dispersion diagram in the $(\omega, k)$ coordinates
($\omega$ is the temporal circular frequency, $k$ is the wavenumber in the $x$-direction).
Usually the dispersion diagram is a graph consisting of several
branches (curves). Each point of the dispersion diagram is characterized by two important parameters, namely
by phase velocity $v_{\rm p} = \omega / k$ and group velocity $v_{\rm g} = (d k / d\omega )^{-1}$ (the derivative
$d k / d\omega$ is the slope of
corresponding branch of the dispersion diagram). It is well known that wave pulses  propagate in the waveguide with corresponding group velocities.

It is well known also that in the experiment typically one can observe the {\em first arriving signal} (FAS)
propagating with the velocity of the fastest medium in the structure of the waveguide. In the case of a
single elastic isotropic medium the velocity of the FAS is close to the velocity of the longitudinal waves.
In some cases the velocity of the FAS is bigger than any of the group velocities provided by the dispersion diagram. The FAS pulse decays with propagation distance unlike the pulses corresponding to usual modal pulses (the latter are called {\em guided waves} in the medical-related literature). Thus, FAS is a transient process in a waveguide. It can negligibly small very far from the source, however at moderate distances it can be an important phenomenon.

The most known applications of FAS are related to medical acoustics (see e.g. \cite{Moilanen,Grondin}). A long bone can be considered as a tubular  waveguide constituted of three media: a thin outer layer of dense (cortical) bone, a sponge bone underneath, and a liquid marrow core in the center. The fastest  wave that can theoretically propagate in such media (taken separately as infinite spaces) is the longitudinal wave in the cortical bone. The sponge bone and the marrow bear much slower waves. However, when a standard analysis of a waveguide is performed (say, by the finite element method) the dispersion diagram contains no branch having group velocity close to the longitudinal cortical velocity.

There exist different approaches to describe FAS. The most comprehensive approach has been proposed by Miklowitz and Randles in \cite{MiklowitzRandles}. This approach considers an analytical continuation of the dispersion diagram.
After passing some branch points, it is possible to find a single branch of the dispersion diagram whose
group velocity and decay correspond to FAS. Another analytical approach is an application of the classical Cagniard--deHoop technique to invert the Fourier transform in the time domain \cite{Miklowitz}.
An approximate approach to FAS is described in \cite{Camus}, where this type of waves is treated as a head wave.
Also FAS can be approximately described as a leaky mode. For this, the slow media composing the waveguide are
declared as elastic half-spaces, while the fastest medium remains to be a layer. Such an approach enables one to compute
(at least in the simplest cases) the velocity and the decay of FAS with an acceptable accuracy.

\begin{figure}[ht]
\centerline{\epsfig{file=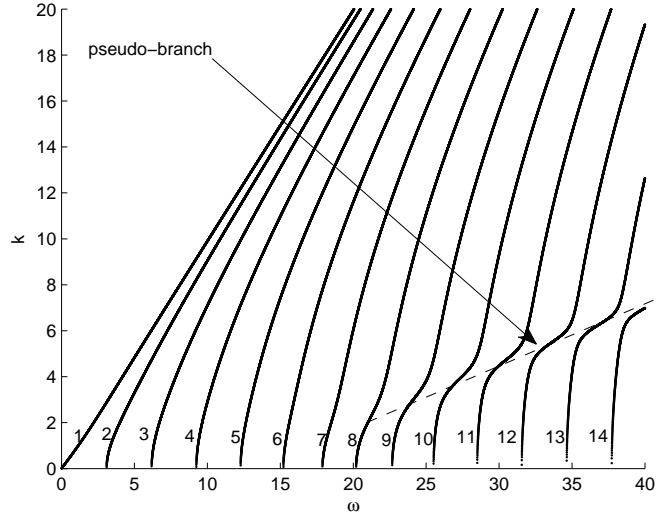,width=10cm}}
\caption{A dispersion diagram of the waveguide composed of two layers}
\label{fig01}
\end{figure}

The aim of the current paper is to present a model of FAS based on an elementary analysis of real dispersion diagrams of a layered waveguide. It is known \cite{Mindlin} that a wave process in a layered waveguide can be treated as
an interaction between the modes of different types and velocities. Thus, the dispersion diagram has a ``terrace-like'' structure formed by overlapping of different sets of branches. Since typically no crossing
of the branches can happen (except the branches corresponding to non-interacting waves), there occur
quasi-crossings at which the type of the mode is changing. A typical fragment of a dispersion diagram is shown in
Fig.~\ref{fig01} (this is a dispersion diagram
of a model two-media acoustic waveguide studied in the paper). The visible line of the smallest slope (although this line can be composed of segments corresponding to different branches) is a {\em pseudo-branch}) relating to the FAS.

In the current paper this idea is developed into an analytical model. An approximation of a fragment of the dispersion diagram by a sum of a tangent function   and a linear function is constructed.
An analytical continuation of the approximation to the domain of positive ${\rm Im}[\omega]$ is performed. The tangent
function tends to an imaginary constant there. Using this trick, the estimation of the decay of FAS is obtained.

The idea to study the analytical continuation of the dispersion diagram is inspired
by the Miklowitz--Randles method \cite{MiklowitzRandles}.
In the paper we briefly describe this method and its connection with the description of the FAS as a leaky wave.

The structure of the paper is as follows. In Section~II a sample problem is formulated. In Section~III a numerical modeling of pulse propagation in the two-media problem is performed. The presence of FAS
and its exponential decay are established.
In Section~IV the approach by Miklowitz--Randles is applied to the waveguide. The branch of the dispersion diagram
responsible for FAS is found.
In Section~V
the properties of FAS are compared to those of a corresponding leaky wave.
In Section~VI the signal in the waveguide is represented in the form of an approximation of the
 phase based on tangent function.
The approximation is analyzed and the parameters of FAS are estimated. The results obtained by Miklowitz--Randles
approach, leaky wave approach and the analysis of the real dispersion diagram are compared.


\section{A sample two--layered waveguide}

Consider a waveguide in the $(x,y)$-plane occupying the strip $-H_1\le y \le H_2$ (see Fig.~\ref{fig03}). The layer
$H_1 \le y \le 0$ is filled with a medium having density and  speed of sound equal to $\rho_1$, $c_1$,
respectively. The layer $ 0 \le y \le H_2$ is filled with a medium with  parameters $\rho_2$, $c_2$.
The wave equations in the media are as follows:
\begin{equation}
c^2_j (\ptl_x^2 + \ptl_y^2) u_j = \ddot{u}_j,
\label{eq0201}
\end{equation}
where $u_j(x,y,t)$, $j = 1,2$ are the field variables (say, acoustical potentials),  notation $\ddot{u}_j$
stands for the second time derivative.

\begin{figure}[ht]
\centerline{\epsfig{file=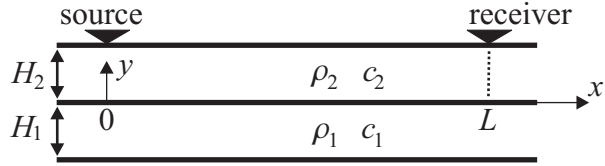}}
\caption{Geometry of a model waveguide}
\label{fig03}
\end{figure}

The boundary conditions are as follows.
The surface $y = -H_1$ is rigid (Neumann):
\begin{equation}
\ptl_y u_1 (x, -H_1, t) = 0.
\label{eq0202}
\end{equation}
On the interface $y = 0$ the pressure and the normal velocity are continuous:
\begin{equation}
\rho_1 u_1 (x,0,t) = \rho_2 u_2 (x,0,t), \qquad \ptl_y u_1 (x,0,t) = \ptl_y u_2 (x, 0, t).
\label{eq0203}
\end{equation}
The surface $y = H_2$ is also rigid, but a point source is located at the point $(0, H_2)$:
\begin{equation}
\ptl_y u_2(x,H_2, t) = \delta(x) f(t),
\label{eq0204}
\end{equation}
where $\delta$ is the Dirac delta-function, $f(t)$ is the time shape of the probe pulse.
The observation point is located at $(L,H_2)$, i.\ e. the function $u_2(L, H_2, t)$ is recorded.

The problem of finding the signal on the receiver is quite standard and it can be solved easily.
Namely, perform the 2D Fourier transform of $u_j$ in the domain of time and frequency:
\begin{equation}
\tilde U_j (k,y,\omega) =
\int \limits_{-\infty}^\infty  \int \limits_{-\infty}^\infty u_j(x,y,t)e^{i\omega t - i k x} dx \, dt.
\label{eq0205}
\end{equation}
Get a 1D problem for $\tilde U_j (k,y,\omega)$ as functions of $y$.
These functions should obey the equations
\begin{equation}
(\ptl_y^2 + \alpha^2_j ) \tilde U_1 (k, y, \omega) = 0,
\qquad
\alpha_j = \alpha_j(k, \omega) = \sqrt{\frac{\omega^2}{c_j^2} - k^2}.
\label{eq0209}
\end{equation}
The following boundary conditions should be valid:
\begin{equation}
\ptl_y \tilde U_1 (k, 0 , \omega) = \ptl_y \tilde U_2 (k, 0 , \omega),
\qquad
\rho_1 \tilde U_1 (k, 0 , \omega) = \rho_2 \tilde U_2 (k, 0 , \omega),
\label{eq0211}
\end{equation}
\begin{equation}
\ptl_y \tilde U_1 (k, -H_1 , \omega) =0 ,
\qquad
\ptl_y \tilde U_2 (k, H_2 , \omega) = F(\omega),
\label{eq0210}
\end{equation}
where
\begin{equation}
F(\omega) = \int \limits_{-\infty}^\infty f(t) e^{i \omega t} dt
\label{eq0205a}
\end{equation}

The solution of (\ref{eq0209}), (\ref{eq0210}), (\ref{eq0211}) can be found:
\begin{equation}
\tilde U_2 (k, H_2, \omega) = F(\omega) \frac{ M(k, \omega)}{N(k, \omega)},
\label{eq0212}
\end{equation}
\begin{equation}
M(k, \omega) = \frac{\alpha_1}{\alpha_2} \sin (\alpha_1 H_1) \sin(\alpha_2 H_2) -
\frac{\rho_1}{\rho_2} \cos(\alpha_1 H_1) \cos(\alpha_2 H_2),
\label{eq0213}
\end{equation}
\begin{equation}
N(k, \omega) = \alpha_2 \frac{\rho_1}{\rho_2} \cos(\alpha_1 H_1) \sin(\alpha_2 H_2) +
\alpha_1 \sin(\alpha_1 H_1) \cos(\alpha_2 H_2).
\label{eq0214}
\end{equation}
The field on the receiver can be obtained by inverting the Fourier transformation:
\begin{equation}
u_2 (x, H_2, t) = \frac{1}{4 \pi^2} \int \! \! \! \! \!  \int \limits_{-\infty}^{\infty}
F(\omega) \frac{ M(k, \omega)}{N(k, \omega)} e^{i k x - i \omega t} d k d \omega
\label{eq0215}
\end{equation}

Formula (\ref{eq0215}) cannot be used directly, since zeros of the denominator belong to the
plane of integration. The limiting absorption principle is used to change the contour of integration.
Namely, for $\omega > 0$ we assume that the velocities $c_j$ have vanishing {\em negative} imaginary parts,
while for $\omega < 0$ the velocities $c_j$ have vanishing {\em positive} imaginary parts. Due to this,
for each $\omega$ the zeros of $N$ in the complex $k$-plane become displaced from the real axis.

The dispersion diagram represents the zeros of $N(k, \omega)$  (for real $c_j$), i.\ e.\ the dispersion
equation is
\begin{equation}
\frac{\alpha_1 \tan(\alpha_1 H_1)}{\alpha_2 \tan(\alpha_2 H_2)} = - \frac{\rho_1}{\rho_2}.
\label{eq0216}
\end{equation}
The roots of (\ref{eq0216}) are curves in the $(k, \omega)$ plane, each point of which
corresponds to a wave freely propagating in the waveguide and
having $x$- and $t$-dependence of the form $~\exp \{ i (kx - \omega t) \}$.


\section{Numerical demonstration of FAS}

The following parameters have been selected
for a numerical demonstration of FAS presence: $H_1 =1$, $H_2 = 0.4$, $c_1 = 1$, $c_2 = 5$,
$\rho_1 = \rho_2 = 1$.
Dimensionless values are used in the computations, since no particular physical medium is under the investigation.
The values plotted in the graphs are, therefore, also dimensionless.
 The dispersion diagram for this waveguide is shown in Fig.~\ref{fig01}.
The fastest pseudo-branch is shown in the figure as a dashed line.  The pseudo-branch is composed of parts of
real branches of the diagram.
One can see that for the selected parameters the pseudo-branch is quite loose, i.~e.\ the
gaps between its parts are quite wide.

For the demonstration we are using the probe pulse $f(t)$ having the spectrum corresponding to the pseudo-branch. Namely, we are using the region $20 < \omega < 40$. In this region the group velocities of the guided waves are
smaller than~$2.7$. In the demonstration we are going to show the presence of a pulse whose velocity is approximately equal to the inverse of the slope
of the dashed line. This velocity is close to  to~4. Thus, the velocity of the FAS is considerably higher than that of
any of the guided waves.

The shape of the pulse $f(t)$
is shown in Fig.~\ref{fig03}, left.
The spectrum of this pulse is shown in Fig.~\ref{fig03}, right.
 One can see that $f(t)$ is a radio pulse centered around $t = 0$. The central circular frequency is
about $\omega_0 = 28$.

\begin{figure}[ht]
\centerline{\epsfig{file=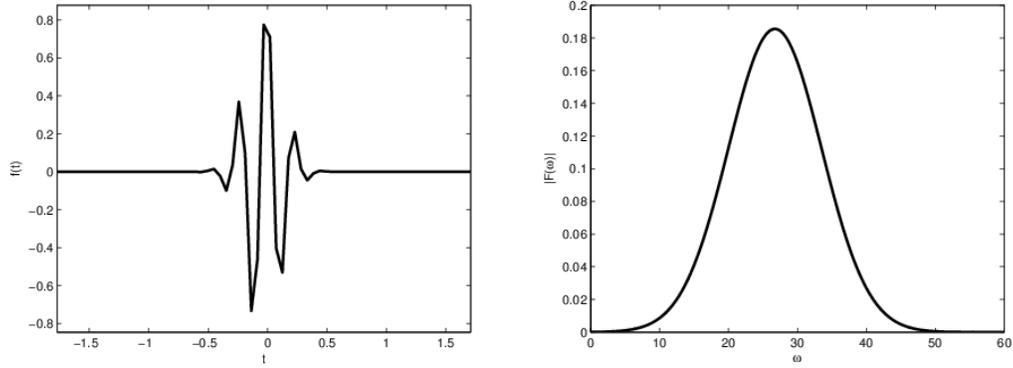,width=14cm}}
\caption{Probing pulse $f(t)$ (left) and its spectrum (right)}
\label{fig03}
\end{figure}

\begin{figure}[h]
\centerline{\epsfig{file=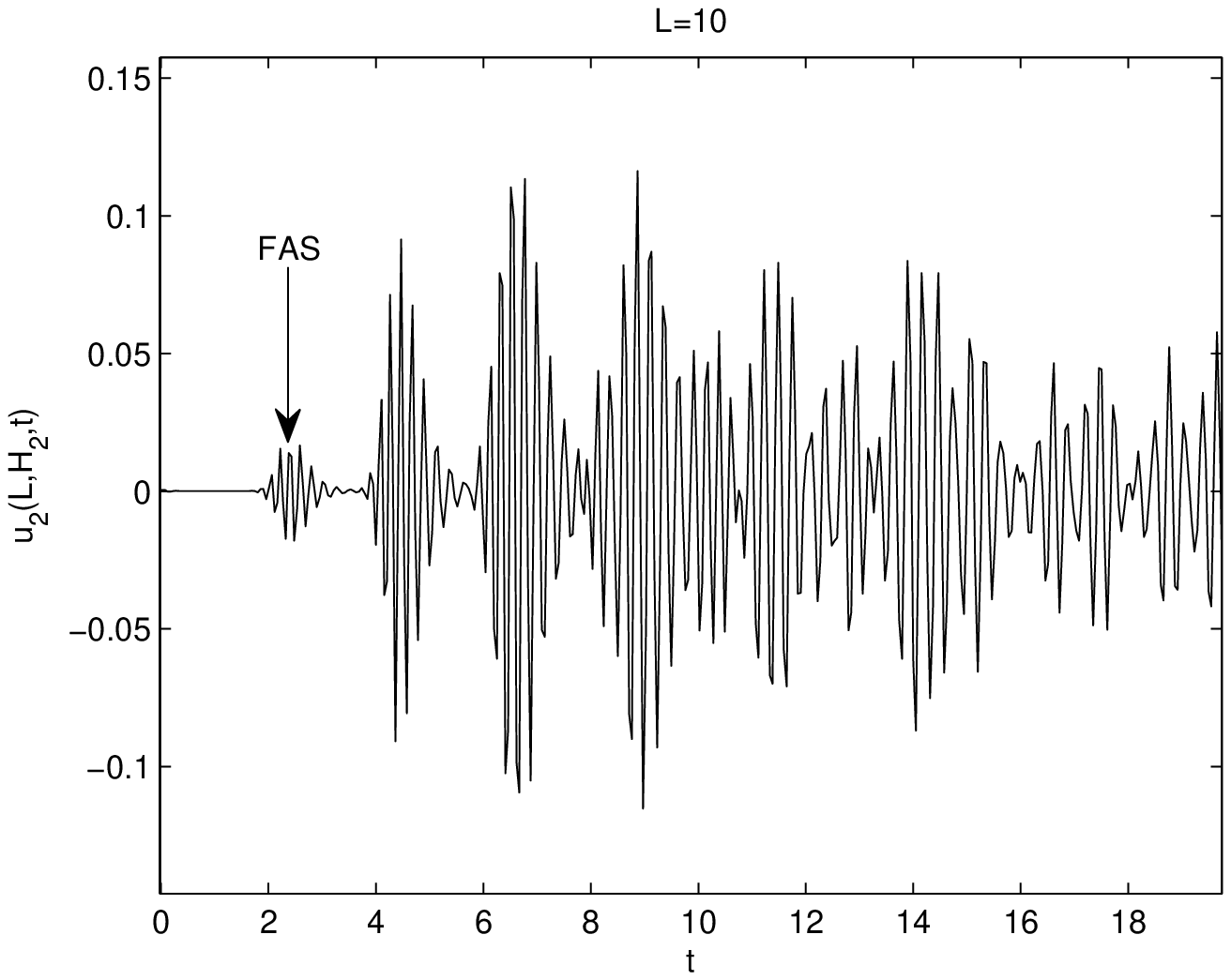,width = 10cm}}
\caption{Output pulse for L = 10}
\label{fig04}
\end{figure}

\begin{figure}[h]
\centerline{\epsfig{file=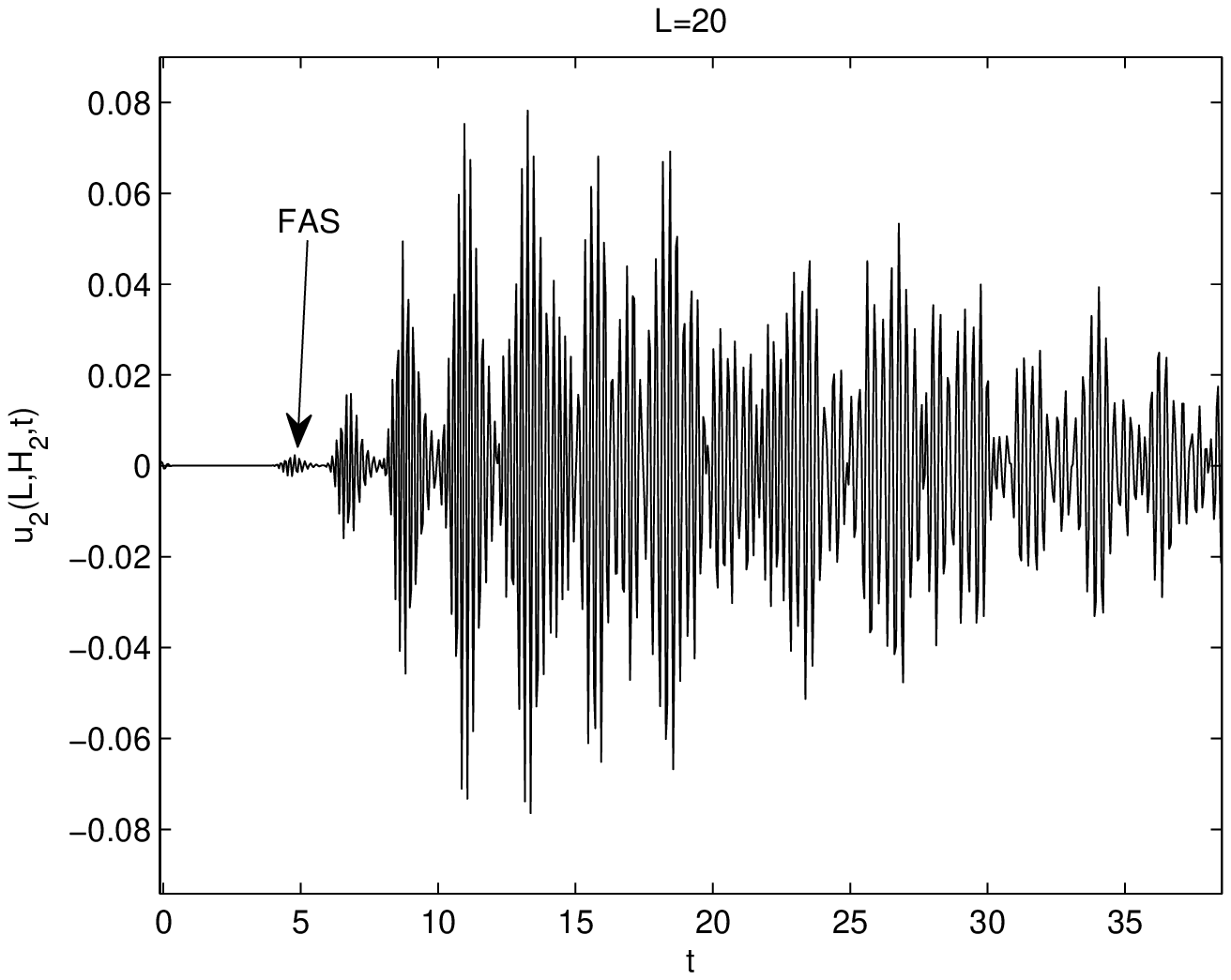,width=10cm}}
\caption{Output pulse for L = 20}
\label{fig05}
\end{figure}

\begin{figure}[h]
\centerline{\epsfig{file=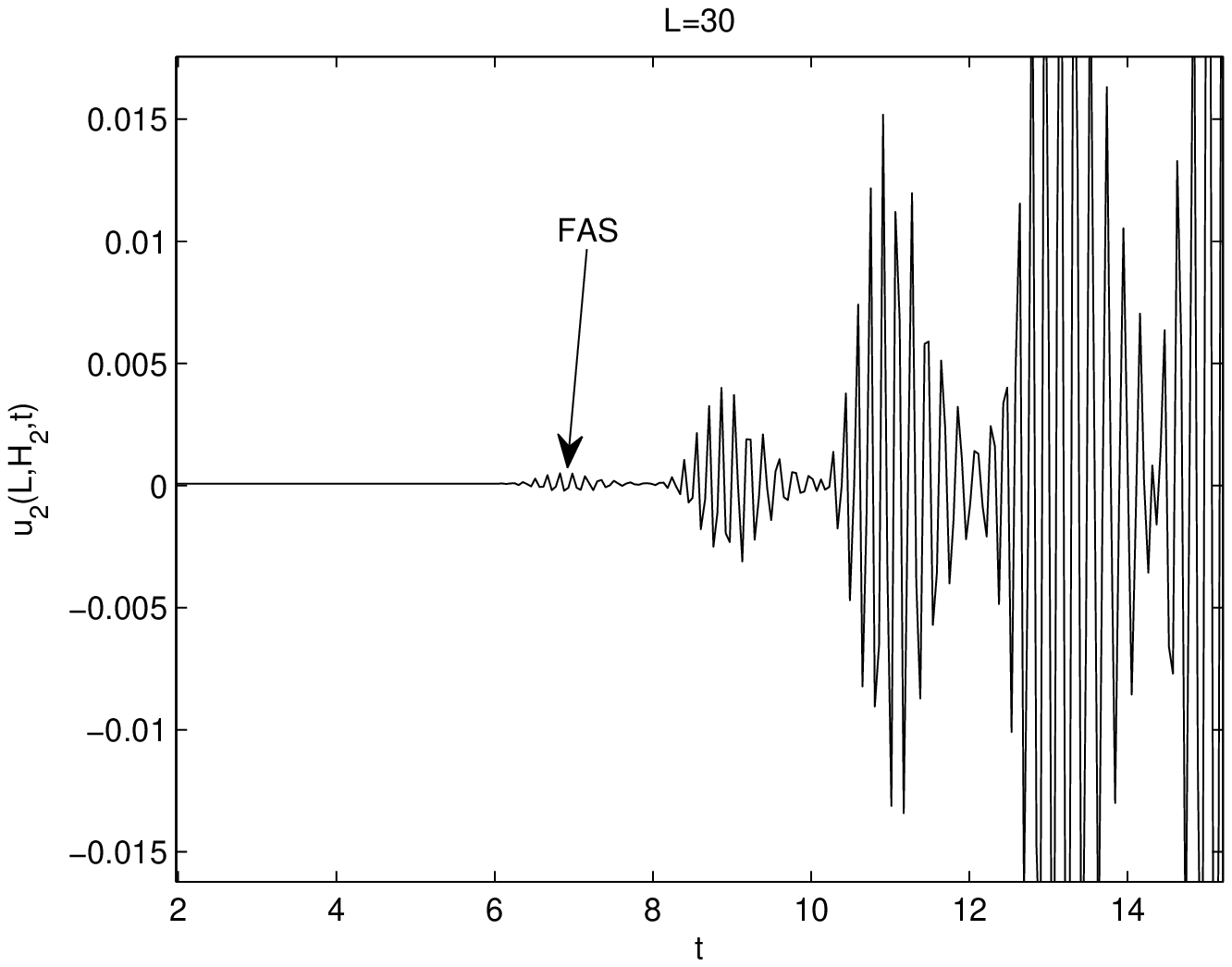,width=10cm}}
\caption{Output pulse for L = 30}
\label{fig06}
\end{figure}

The results of the computations made by formula (\ref{eq0215}) for $L=10,20,30$
are shown in Fig.~\ref{fig04}, Fig.~\ref{fig05}, Fig.~\ref{fig06}, respectively.
The field at the receiver, i.~e. $u_2(L,H_2,t)$, is plotted.
One can see that in all graphs there exists a small pulse, which can be interpreted as FAS.

Parameters of FAS approximately determined from these graphs are given in Table~\ref{tab01}. ToF is the ``time of flight'', i.\ e.\ the travel time of FAS. One can see that the velocity of the pulse is more than 4, and the amplitude decay is close to exponential. Such a behavior is typical for FAS. The aim of the rest of the paper will be to estimate the group velocity of the FAS pulse $v_{\rm g,FAS}$ and the decay parameter $\kappa$. The decay parameter is the coefficient in the exponential dependence of the amplitude vs $L$:
\[
{\rm Amplitude}\sim e^{-\kappa L}.
\]

\begin{table}[ht]
\begin{center}
\begin{tabular}{|c|c|c|}
\hline
$L$ & ToF & Amplitude \\
\hline
10  & 2.5 & $1.8\cdot 10^{-2}$ \\
20  & 4.8 & $2.5\cdot 10^{-3}$ \\
30  & 6.9 & $5.0\cdot 10^{-4}$ \\
\hline
\end{tabular}
\end{center}
\caption{Parameters of the FAS}
\label{tab01}
\end{table}


\section{An approach by Miklowitz and Randles}

Here we are describing a modified version of the approach introduced by Randles and Miklowitz in \cite{MiklowitzRandles}.
The modifications (compared with \cite{MiklowitzRandles}) are as follows: the change of variables made by the authors is omitted, and the variation of the
contour of integration is performed locally in the region of interest with respect to temporal and spatial frequencies.

Let function $f(t)$ be real. Consider the function
\begin{equation}
u_2' (x, H_2, t) = \frac{1}{4 \pi^2} \int\limits_{0}^{\infty}   \int \limits_{-\infty}^{\infty}
F(\omega) \frac{ M(k, \omega)}{N(k, \omega)} e^{i k x - i \omega t} d k d \omega,
\label{eq0401}
\end{equation}
i.\ e.\  exclude the negative values of $\omega$. Obviously,
\begin{equation}
u_2 (x, H_2 , t) = 2{\rm Re}[u_2'(x, H_2, t)].
\label{eq0402}
\end{equation}

For each fixed positive $\omega$ find the roots of (\ref{eq0216}) by solving it as an
equation with respect to~$k$. Denote by $\xi_n (\omega)$ the roots having positive imaginary part or
zero imaginary and positive real part. These roots correspond to waveguide modes traveling in the positive
$x$-direction or decaying in this direction. Note that for for each $\omega > 0$ and for $x > 0$ the integral
with respect to $k$ in (\ref{eq0401}) can be considered as a contour integral, and the contour (the real axis)
can be closed in the upper half-plane. The integrand is a meromorphic function in the upper half-plane,
so the integral can be converted into a sum of residual terms:
\begin{equation}
u_2' (x, H_2, t) = \frac{i}{2 \pi} \int\limits_{0}^{\infty}  F(\omega)  \sum_n
 \frac{ M(\xi_n(\omega), \omega)}{N'(\xi_n (\omega), \omega)} \exp\{i \xi_n(\omega) x - i \omega t\}  d \omega,
\label{eq0403}
\end{equation}
where
\begin{equation}
N' (k, \omega) = \ptl_k N(k, \omega).
\label{eq0404}
\end{equation}
Indeed, (\ref{eq0403}) is a standard expansion of the wave field in the waveguide as a sum of waveguide
modes.

For large $t$ and $x$ (i.\ e.\ in the far field) expansion (\ref{eq0403}) provides a comprehensive
description of the wave field. Namely, each branch of the dispersion diagram with real
$\xi$ corresponds to a modal pulse, and the
velocities of the modal pulses are provided by the group velocities of the branches.
The branches with imaginary $\xi$ form the near field.

For moderate $t$ and $x$, however, some transient processes can occur in the waveguide. The most
interesting example of such transient processes is the FAS demonstrated in the previous section.
In the transient region the standard analysis of the real dispersion diagram cannot be performed for two reasons:
\\
\noindent
a) The modal pulses correspond to rather long fragments of the branches, on which one cannot approximate
the branch by taking just its slope and curvature.
\\
\noindent
b) Numerous branches participate in formation of a single pulse.

It was a brilliant idea of Miklowitz and Randles that one can simplify the consideration by deforming the contour of integration of (\ref{eq0403})  into the complex domain of $\omega$. After passing some branch points the structure of the (complexified) dispersion diagram become simpler, and each transient process becomes attributed to a single branch of it.
Let us illustrate this approach on the example of the FAS. As we assume, the FAS demonstrated above corresponds to a fragment of the pseudo-branch shown in Fig.~\ref{fig01}.

Let $\omega$ be a complex variable. Continue functions $\xi_n$ from the real axis to the complex plane.
Consider $\xi_n (\omega)$ as branches of a single analytical multivalued function $\xi (\omega)$. This function, obviously, has an infinite number of sheets. Its analyticity follows
from standard theorems applied to equation (\ref{eq0216}). Function $\xi(\omega)$ is defined on its Riemann surface.
Practically, the continuation of $\xi_n$ is obtained by numerical solving of equation (\ref{eq0216}) with respect to
$k$ for complex $\omega$.

\begin{figure}[ht]
\centerline{\epsfig{file=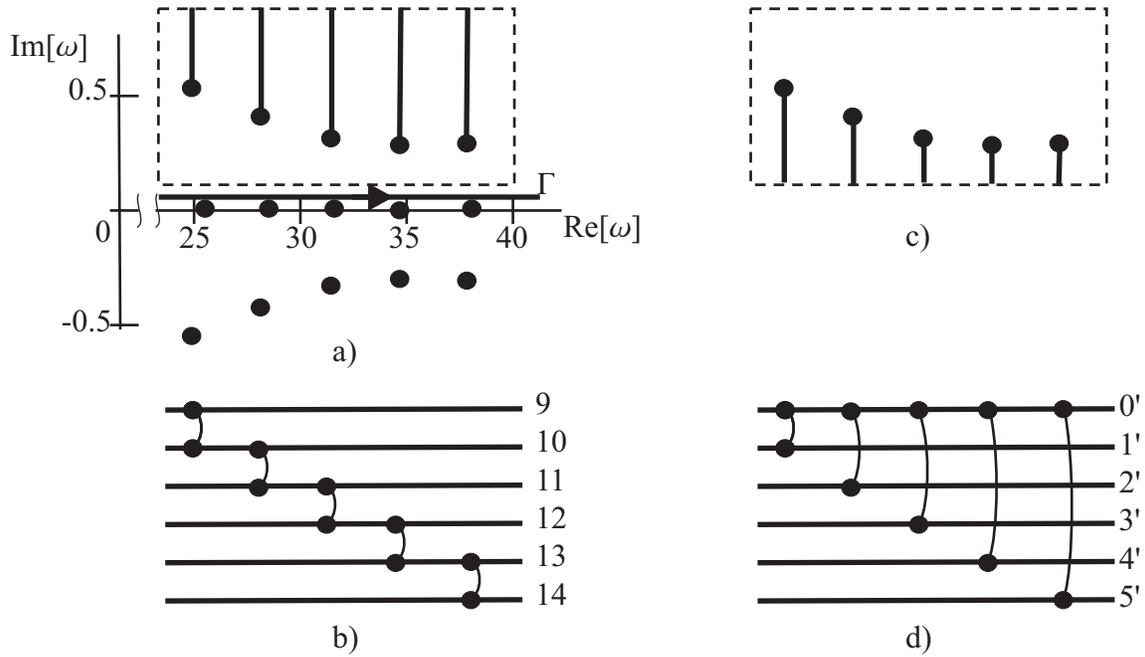}}
\caption{Branch points on a fragment of the Riemann surface of $\xi(\omega)$}
\label{figA01}
\end{figure}

Explore a fragment of the Riemann surface. The area under investigation is shown in Fig.~\ref{figA01}a. We are interested in behavior of the branches corresponding to the real curves labeled by numbers 9--13 in Fig.~\ref{fig01}.
The fragment of the Riemann surface is shown in Fig.~\ref{figA01}a.
Small circles denote the branch points. All branch points are of order~2, i.~e.\ they connect two sheets of the surface. The branch points form three groups: the points belonging to the real axis, the points above the
real axis, and the points below the real axis.

The points on the real axis correspond to the cut-off points of the waveguide modes, i.~e.\ these branch points
connect propagating branches with evanescent branches.

The branch points above the real axis play an important role. They transform a pseudo--branch into a leaky wave branch.
Consider the area of the Riemann surface bounded by a dashed rectangle in Fig~\ref{figA01}a. Cut the fragments of the sheets lying in this area  by drawing the
cuts from the branch points to the top side of the rectangle (the cuts are shown by bold lines). The scheme of the Riemann surface is shown in Fig~\ref{figA01}b. The indices of the sheets correspond to the branches of the dispersion diagram in Fig~\ref{fig01}.

 This scheme reflects the structure of the pseudo-branch, namely, the sheets are linked sequentially.
Each sheet is linked with the previous one and with the next one.
The same surface can be cut in a different way. Make the cuts going from the branch points to the bottom side of the rectangle. (see Fig.~\ref{figA01}c). As the result, get the surface having scheme shown in Fig.~\ref{figA01}d. One can see that
there is a single sheet (labeled by $0'$), to which all other sheets are linked. This sheet  bears the
branch of the dispersion diagram similar to the leaky wave branch described in the next section.

\begin{figure}[ht]
\centerline{\epsfig{file=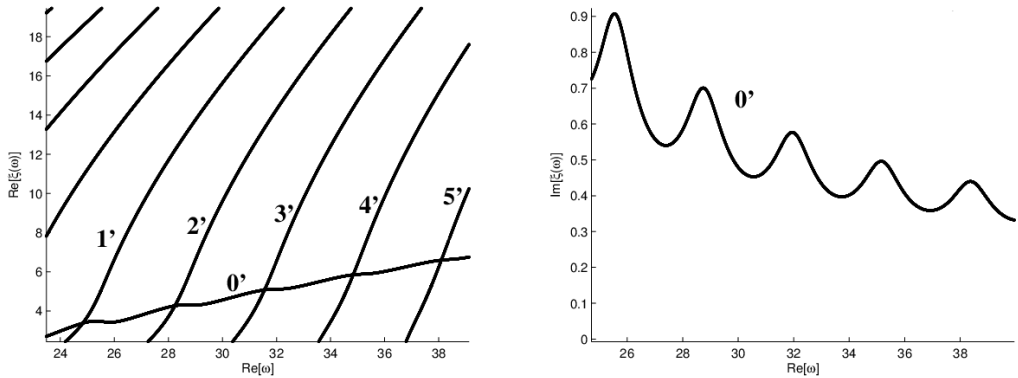,width=14cm}}
\caption{A fragment of the analytical continuation of the dispersion diagram
$\xi(\omega)$ for ${\rm Im}[\omega]=1$}
\label{figA03}
\end{figure}

To illustrate the structure of the Riemann surface shown in Fig.~\ref{figA01}d, plot several branches of $\xi(\omega)$
for $\omega$ having positive imaginary part higher than the branch points. In Fig.~\ref{figA03}
(left) we plot the real part of the dispersion diagram for ${\rm Im}[\omega] = 1$. The indices
correspond to the notations of sheets
in Fig~\ref{figA01}d.
One can see that the branch $0'$ has a small
slope $d\xi/d\omega$ i.~e.\ a high group velocity, and there are some other branches with smaller group velocity. In the right part of the figure we plot the imaginary part of $\xi(\omega)$ corresponding to the sheet~$0'$.
The imaginary part of $\xi$ for other branches are much higher, i.~e.\ the corresponding wave components have
stronger attenuation.

Let us use the knowledge of the Riemann surface of $\xi(\omega)$ to analyze integral
(\ref{eq0403}). Rewrite this integral in the form
\begin{equation}
u_2' (x, H_2, t) = \frac{i}{2 \pi} \int_\Gamma A(\omega) \exp\{i \xi(\omega) x - i \omega t\}  d \omega,
\label{eq0403a}
\end{equation}
where contour $\Gamma$ is a sum of an infinite number of contours going from $0$ to $\infty$ along different sheets of the
Riemann surface of $\xi(\omega)$, and $A(\omega)$ is a multivalued function (all non-exponential
factors of the integrand of (\ref{eq0403})). A natural assumption can be made that $A(\omega)$ can be continued into the upper half-plane
of $\omega$, and the Riemann surface of $A(\omega)$ has the same branch points and the same topology as the Riemann surface of $\xi(\omega)$.

\begin{figure}[ht]
\centerline{\epsfig{file=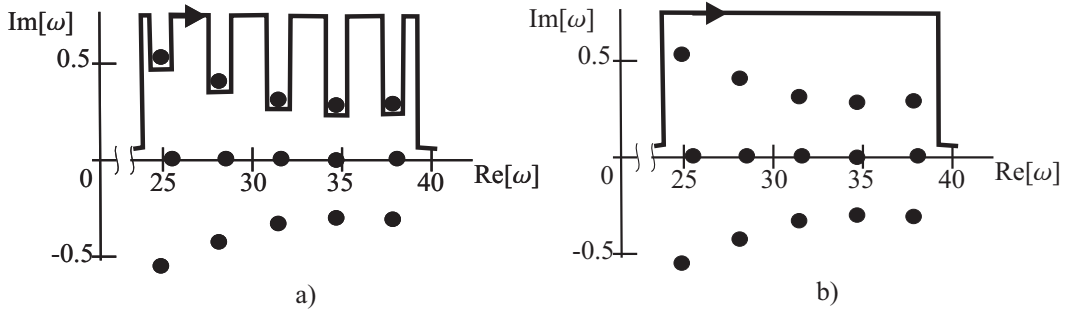,width=14cm}}
\caption{Deformation of contour of integration in (\ref{eq0403a})}
\label{figA02}
\end{figure}

The idea of Miklowitz and Randles (formulated in a slightly different form and for a different type of waveguide)
is to deform the contour of integration $\Gamma$ shown in Fig.~\ref{figA01}a first into the contour
shown in Fig.~\ref{figA02}a and then into the contour shown in Fig.~\ref{figA02}b.
Note that $\Gamma$ is a set of contours on different sheets, and all of them are deformed simultaneously.

The change of contour of integration should be commented as follows.

\begin{itemize}

\item
Transformation from the contour shown in Fig.~\ref{figA02}a to the contour shown in
Fig.~\ref{figA02}b
is possible
since all branch points have order~2, and contour $\Gamma$ is composed of contours
located on all sheets. The integrals along the vertical parts encircling the branch points compensate each other.

\item
For some $x$ and $t$ such a deformation leads to an exponential decrease of the exponential factor of the integrand.
Namely, let $x/t$ be greater than any of the values $({\rm Re}[d \xi /d\omega ])^{-1}$ within the area in the dashed
rectangle. These values are indeed the group velocities of the branches. Then the exponential factor
$\exp \{ i x (\xi  - i \omega t/x) \}$ decreases as ${\rm Im}[\omega]$ grows, and this decrease is exponential
if $x$ is a large parameter for fixed $x/t$. Thus, the contour deformation can be used for finding the transient
wave components that are faster than usual modal pulses.

\end{itemize}

After the deformation of the contour one gets the integral on the sheet labeled as $0'$ (describing the FAS)
and many other
integrals corresponding to smaller wave components, which can be neglected. Since the behavior of the branch $0'$ is more
regular than that of the initial real dispersion diagram, one can introduce the group velocity of FAS by
\begin{equation}
v_{\rm g,FAS} = \left( \frac{d \xi_{0'}}{ d\omega}  \right)^{-1}
\label{eq0403b}
\end{equation}
and define the decay parameter $\kappa$ of the FAS
as
\begin{equation}
\kappa = {\rm Im}[\xi_{0'}(\omega)] - \frac{{\rm Im}[\omega]}{v_{\rm g,FAS}}
\label{eq0404}
\end{equation}
taken for ${\rm Re}[\omega] = \omega_0$, which is medium for the wave train. Anyway, an exact definition of group velocity and decay is impossible, since generally FAS is a highly dispersive pulse.

The graphs shown in Fig.~\ref{figA03} enable one to estimate $v_{\rm g,FAS}$ as approximately $4.0$.
Parameter $\kappa$ varies considerably within the frequency band of the pulse (between $0.1$ and $0.6$).
One can assume that the components of the $FAS$ with higher frequencies has lower decay than the low-frequency
components.


\section{FAS as a leaky wave}

\begin{figure}[ht]
\centerline{\epsfig{file=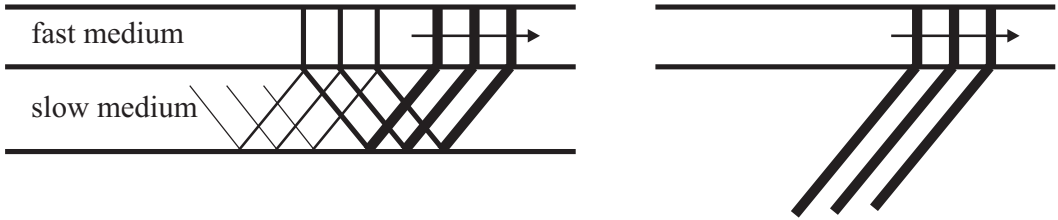,width=14cm}}
\caption{Structure of FAS (left) and leaky wave (right)}
\label{figA04}
\end{figure}

To study the leaky wave we make the lower medium occupying the whole half-space $y<0$.
Let us explain why the analytical continuation of the dispersion diagram into the complex domain of $\omega$ reveals
a mode whose structure is close to a leaky wave. The structure of FAS in a simplest case is shown in Fig.~\ref{figA04}. One can see that it consists of a leading pulse traveling in the fast medium and the head wave traveling in the slow medium.
The dispersive relation for FAS is similar to that of the leaky wave if the slow medium can be
substituted by a half-space, i.\ e.\ if the waveguide does not feel the lower boundary.
If the analytical continuation into the complex domain of $\omega$ is made, and if $k \approx \omega / c_2$ (which is the case), then $\alpha_0 \approx 0$, while $\alpha_1$ has a considerable imaginary part. The imaginary part of $\alpha_1$
corresponds to transversal decay of the waves in the slow medium. This decay eliminates the influence of the lower boundary.

The structures
of the FAS and the leaky wave
shown in Fig.~\ref{figA04} demonstrates that the leaky wave model should not work well for strongly
dispersive FAS. Namely, if the FAS pulse in the upper medium starts to be wider than the distance between the beginning of the pulse and the reflection coming from the lower boundary, the leaky wave model fails to describe the wave process.

\begin{figure}[ht]
\centerline{\epsfig{file=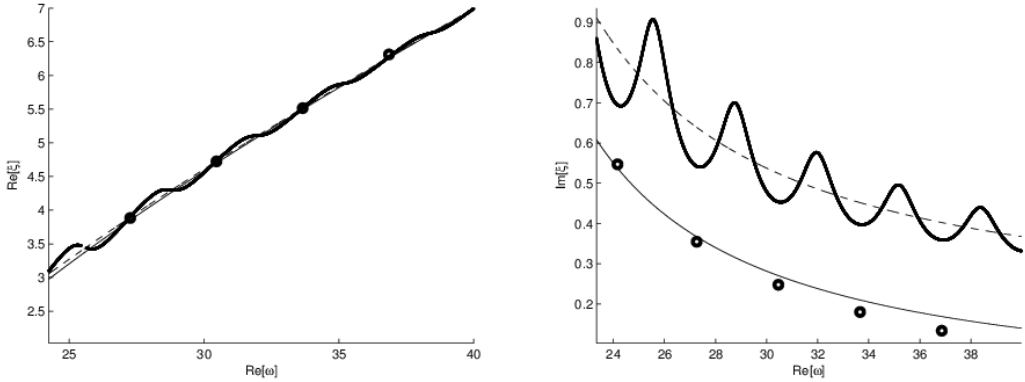,width=14cm}}
\caption{Comparison of the leaky wave dispersion diagram with the branch $\xi_{0'}$ and with the estimations based on the
real dispersion curve. The left side displays
the real part, the right side displays the imaginary part}
\label{figA05}
\end{figure}

The dispersion relation for the leaky wave in our case has form
\begin{equation}
\frac{\rho_1 \alpha_2(\omega, \xi_{\rm L}) \tan(\alpha_2(\omega , \xi_{\rm L}) H_2)}
{\rho_2 \alpha_1(\omega, \xi_{\rm L})}  = -i
\label{eq0403c}
\end{equation}
In Fig.~\ref{figA05}
the roots of this dispersion relation are compared with analytical continuation of the dispersion diagram constructed before. Thin solid lines correspond to the leaky wave dispersion diagram $\xi_{\rm L}(\omega)$ for
${\rm Im}[\omega] = 0$. Dashed lines correspond to the analytical continuation of the leaky wave dispersion diagram,
namely to ${\rm Im}[\omega] = 1$. The bold lines correspond to the branch $\xi_{0'} (\omega)$ found above for
${\rm Im}[\omega]=1$. The circles are related to the real dispersion diagram analysis described in the next section.

The dispersion relation $\xi_{\rm L} (\omega)$ for the leaky wave can be used for estimation of the parameters of FAS
in a very straightforward way. Namely, the velocity of FAS can be estimated as ${\rm Re}[d\xi_{\rm L} / d\omega]^{-1}$,
while the decay parameter $\kappa$ can be approximated as $\kappa \approx {\rm In}[\xi_{\rm L} (\omega_0)]$.

\section{Analysis of the real dispersion diagram}

In the previous sections we described two approaches to FAS. The Miklowitz--Randles approach requires a continuation of the dispersion diagram into the domain of complex $\omega$. The leaky wave approach requires a reduced physical model
that can be sometimes difficult to construct. Here we describe an alternative approach based on analysis of the real dispersion diagram. One of the benefits of this approach is that it can be used even when the only available information
about the waveguide is a numerical description of the dispersion diagram.
However, the concepts described above are necessary for understanding of the current approach.

Consider an approximation to $\xi(\omega)$ taken in the vicinity of the pseudo-branch:
\begin{equation}
\xi_{\rm a} (\omega) = a \tan (\beta \omega + \gamma) + \frac{\omega}{v} + c
\label{eq1001}
\end{equation}
where $a$, $\beta$, $\gamma$, $v$, $c$ are some parameters.
Parameter $\beta$ is linked to the number of real branches crossing the pseudo-branch per some unit length
$\Delta \omega$. Parameter $a$ shows how loose is the pseudo-branch.

Consider the field in vicinity of the
point $(t, x = v t)$. To study the integral of the form (\ref{eq0403a}) deform the contour of integration
as it was described above and continue $\xi(\omega)$ into the upper half-plane of complex $\omega$.
The most important observation is that for large $\beta$
\[
\tan(\beta \omega + \gamma) \to i
\]
in the upper half-plane.
Thus, the decay of the wave corresponding to the pseudo-branch can be estimated as
\begin{equation}
\kappa \approx a.
\label{eq1002}
\end{equation}
Obviously, the velocity of the wave component associated with the pseudo-branch is equal to $v$.

The procedure of estimation of $v$ and $a$ is as follows. It is natural to assume that the pseudo-branch
passes through the points of the dispersion diagram, at which the slope $d\xi / d\omega$ has local minimums
(i.~e.\ the group velocity has local maximums). Denote these points by $\omega^*_n$.
We are trying to find such parameters $v$ and $c$ that the linear function $\xi = \omega/v + c$ fits
these points of the dispersion diagrams. Then we try find the tangent function parameters $a$, $\beta$ and $\gamma$
such that the approximation function (\ref{eq1001}) coincides with the real dispersion diagram $\xi(\omega)$
at the points $(\omega^*_n , \xi(\omega^*_n))$ and has the slope $d \xi_{\rm a}(\omega^*_n) /d \omega$
at these points equal to the slope $d \xi(\omega^*_n) /d \omega$ of the dispersion diagram.

\begin{figure}[ht]
\centerline{\epsfig{file=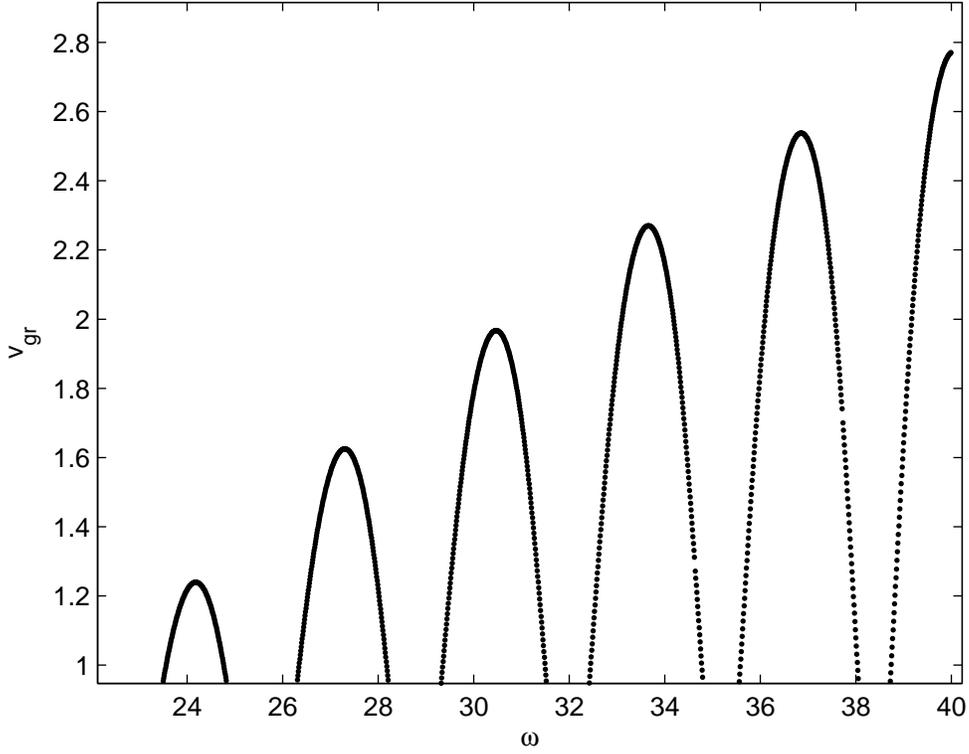}}
\caption{Group velocities of the modes}
\label{figA06}
\end{figure}

In Fig.~\ref{figA06} we plot a fragment the group velocities
graph related to Fig.~\ref{fig01}. The wavenumber $k$ and the group velocity $v_{\rm gr}$ related to each
local maximum $\omega^*_n$ are given by Table~\ref{tab02}.

\begin{table}[ht]
\begin{center}
\begin{tabular}{|c|c|c|c|c|c|}
\hline
$\omega_n^*$                  & 24.16 & 27.26 & 30.46 & 33.66 & 36.86  \\
\hline
$k(\omega_n^*)$               & 3.0   & 3.88  & 4.72  & 5.51  & 6.31   \\
\hline
$v_{\rm gr}(\omega_n^*)$      & 1.24  & 1.63  & 1.97  & 2.27  & 2.54   \\
\hline
\end{tabular}
\end{center}
\caption{Characteristics of local peaks of group velocity}
\label{tab02}
\end{table}

The estimation of $v$ is, obviously,
\begin{equation}
v \approx \frac{\omega^*_{n+1} - \omega^*_n}{k(\omega^*_{n+1}) - k(\omega^*_{n+1})}
\label{eq1003}
\end{equation}
This parameter is an estimation of $v_{\rm g, FAS}$.
The estimation of $\beta$ is
\begin{equation}
\beta \approx \frac{\pi}{ \omega_{n+1}^* -  \omega_{n}^* }.
\label{eq1004}
\end{equation}
Finally, the estimation of $a$ is as follows:
\begin{equation}
a \approx \frac{1}{\beta} \left(
\frac{d \xi (\omega^*_n)}{d \omega} - \frac{1}{v}
\right) .
\label{eq1005}
\end{equation}
Since the values of $a$ given by (\ref{eq1005}) are different for different points $\omega^*_n$, we come to conclusion that $a$ in (\ref{eq1001}) is a (slow) function of $\omega$. By (\ref{eq1005}) we obtain an estimation of $a$ at the
points $\omega_n^*$, and after that we can, say, interpolate $a$ between these points.
As we mentioned above, function $a$ is an estimation of the decay parameter $\kappa$.

Using  the data extracted from Fig.~\ref{figA06} and Fig.~\ref{fig01}
collected in Table~\ref{tab02},
we can estimate the position
of the pseudo-branch and the values of the attenuation parameters. The points related to the peaks of the
group velocity are plotted in Fig~\ref{figA05} as small circles. One can see that the points obtained from the the real dispersion diagram
are in reasonably good agreement with the dispersion curve of the leaky wave for real $\omega$.


\section{Conclusion}

It is shown that FAS pulses correspond to terrace-like structures of the dispersion diagrams. These structures
are called pseudo-branches in the paper. A typical form of the pseudo-branch is shown in Fig.~\ref{fig01}.
The pulse propagation velocity for FAS is the inverse of the slope of the dashed curve in Fig.~\ref{fig01}, and it can be larger than any of the group velocities available in the considered part of the dispersion diagram.

In the paper we describe three different approaches to finding the velocity and the decay of FAS pulses. They are
the Miklowitz--Randles approach based on the analytical continuation of the dispersion diagram, The leaky wave approach
and the approach based on the analysis of the pseudo-branch of the real dispersion diagram. Only the third approach is
original. In the third approach one should study the local peaks of the the group velocity, i.~e.\ find the positions of the peaks and the values of the group velocity. These parameters provide data describing the FAS.
Our computations prove consistency of all three approaches.

Author is grateful to an anonymous referee who attracted his attention to the leaky wave approach.

The work is supported by the grants RFBR 14-02-00573 and Scientific Schools-283.2014.2.



\end{document}